# Oscillation and stability of first-order delay differential equations with retarded impulses

Başak KARPUZ[*]

**Address**. Department of Mathematics, Faculty of Science and Literature, ANS Campus,
Afyon Kocatepe University, 03200 Afyonkarahisar, Turkey

## Abstract

In this paper, we study both the oscillation and the stability of impulsive differential equations when not only the continuous argument but also the impulse condition involves delay. The results obtained in the present paper improve and generalize the main results of the key references in this subject. An illustrative example is also provided.

## 1 Introduction

The theory of impulsive differential equations is an important area of scientific activity, since every nonimpulsive differential equation can be regarded as an impulsive differential equation with no impulse effect, i.e., the corresponding impulse factor is the unit. This fact makes more interesting than the corresponding theory of nonimpulsive differential equations. Moreover, such equations naturally appear in the modeling of several real-world phenomena in many areas such as physics, biology and engineering.

First paper on oscillation of impulsive delay differential equations [1] was published in 1989. From the publication of this paper up to the present time, impulsive delay differential equations started receiving attention of many mathematicians and numerous papers have been published on this class of equations. Most of the publications are devoted to oscillation of first-order impulsive delay differential equations with instantaneous impulse conditions (see for instance [1, 2, 3, 4, 5, 6, 7, 8, 9, 10]). However, to the best of our knowledge, there is not much done in the direction of oscillation and stability of impulsive delay differential equations when the impulse condition also involves a delay argument. Results dealing with retarded impulse conditions are relatively scarce, for instance, we can only a few papers which only deal with the stability property of the solutions (see for instance [11, 12]). Unfortunately, there seems to be nothing accomplished in revealing the oscillation properties of equations with retarded impulse conditions. In this paper, we shall draw our attention to this untouched problem. More precisely, the aim of this paper is to deliver answers to the following questions on the stability and oscillatory behaviour of the solutions to impulsive delay differential equations: "If delay differential equations without impulses are oscillatory, will their solutions continue to oscillate in the absence of retarded impulse perturbations?" and "If the zero solution of delay differential equations without impulses is stable, when exposed to retarded impulse effects, under what conditions the equation maintains the stability?" The motivation of this paper mainly originates from the work [7], where Yan and Zhao studied impulsive delay differential equations of the form

$$\begin{cases} x'(t) + p(t)x(\tau(t)) = 0 & \text{for } t \in [\theta_0, \infty) \backslash \{\theta_k\}_{k \in \mathbb{N}_0} \\ x(\theta_k^+) = \lambda_k x(\theta_k) & \text{for } k \in \mathbb{N}_0 \end{cases}$$

and established very important connections with the following nonimpulsive delay differential equation

$$y'(t) + \left[\prod_{\tau(t) \leq \theta_k < t} \frac{1}{\lambda_k}\right] p(t) y(\tau(t)) = 0 \quad \text{for } t \in [\theta_0, \infty) \text{ almost everywhere.}$$

---

[*]**Email**: bkarpuz@gmail.com **Web**. http://www2.aku.edu.tr/~bkarpuz

Practically, their method introduces a transform which glues the continuous pieces in the graph of a jump type discontinuous function endwise, i.e., sticks together the points at which the jump magnitude of the function is momentary. The method developed by the authors is very effective since it allows handling both oscillation and stability of impulsive differential equations. But, due to technical and theoretical obstacles, the same method is useless for studying impulsive delay differential equations when the impulse condition involves delays as well. We shall therefore introduce a new method which generalizes the method in [7] to such problems. Roughly speaking, the technique still relies on construction of continuous functions from a jump type discontinuous function but unlike the previous case the jump magnitude at each given time is now allowed to depend on the magnitude of a prior time. The technique employed in the present paper allows us to study of both oscillation and stability of such equations without putting any sign condition on the coefficient. In [13, 14, 15], the readers may find stability results for delay differential equations without impulses, which can be combined with our results. The results of this paper improve, generalize and extend the qualitative theory of differential equations to impulsive differential equations with retarded impulse conditions. We refer the readers to the book [16], which covers the fundamental results of the theory on impulsive differential equations.

Our attention in this paper centers on the qualitative behaviour of solutions of the impulsive delay differential equation

$$\begin{cases} x'(t) + p(t)x(\tau(t)) = 0 & \text{for } t \in [\theta_0, \infty) \backslash \{\theta_k\}_{k \in \mathbb{N}_0} \\ x(\theta_k) = \lambda_k x(\theta_{k-\ell}^-) & \text{for } k \in \mathbb{N}_0 \end{cases} \quad (1)$$

under the following primary assumptions:

(A1) $p : [\theta_0, \infty) \to \mathbb{R}$ is a Lebesgue measurable and a locally essentially bounded function.

(A2) $\tau : [\theta_0, \infty) \to \mathbb{R}$ is a Lebesgue measurable function satisfying $\tau(t) \leq t$ for all $t \in [\theta_0, \infty)$ and $\lim_{t \to \infty} \tau(t) = \infty$.

(A3) $\ell \in \mathbb{N}_0$ and $\{\theta_k\}_{k \in \mathbb{N}} \cup \{\theta_{-\ell}, \theta_{-\ell+1}, \ldots, \theta_0\}$ is an increasing divergent sequence of reals.

(A4) $\{\lambda_k\}_{k \in \mathbb{N}_0}$ is a sequence of reals which has no zero terms.

We define $AC([\rho_{\theta_0}, \theta_0], \mathbb{R})$, where $\rho_t := \inf\{\tau(\xi) : \xi \in [t, \infty)\}$ for $t \in [\theta_0, \infty)$, to be the set of functions absolutely continuous functions defined on $[\rho_{\theta_0}, \theta_0]$.

**Definition 1** (Solution). *Suppose that (A1)–(A4) hold. A function $x : [\rho_{\theta_0}, \infty) \to \mathbb{R}$ denoted by $x(\cdot, \theta_0, \varphi)$ is called a* solution *of the initial value problem*

$$\begin{cases} x'(t) + p(t)x(\tau(t)) = 0 & \text{for } t \in [\theta_0, \infty) \\ x(\theta_k) = \lambda_k x(\theta_{k-\ell}^-) & \text{for } k \in \mathbb{N}_0 \\ x(t) = \varphi(t) & \text{for } t \in [\rho_{\theta_0}, \theta_0), \end{cases} \quad (2)$$

*where $\varphi \in AC([\rho_{\theta_0}, \theta_0], \mathbb{R})$ is given, if the following conditions are satisfied:*

(i) *For any $k \in \mathbb{N}_0$, $x$ is absolutely continuous on the interval $[\theta_k, \theta_{k+1})$.*

(ii) *For any $k \in \mathbb{N}_0$, both right-sided and left-sided limits of $x$ exist at $\theta_k$ with $x(\theta_k) = x(\theta_k^+)$.*

(iii) *$x$ is equal to the initial function $\varphi$ on the interval $[\rho_{\theta_0}, \theta_0]$, and satisfies the differential equation*

$$x'(t) + p(t)x(\tau(t)) = 0 \quad \text{for every } t \in [\theta_0, \infty) \setminus \{\theta_k\}_{k \in \mathbb{N}_0}.$$

(iv) *$x$ satisfies the impulse condition*

$$x(\theta_k) = \lambda_k x(\theta_{k-\ell}^-) \quad \text{for all } k \in \mathbb{N}_0,$$

*and it may have jump type discontinuity at the impulse points $\{\theta_k\}_{k \in \mathbb{N}_0}$.*

**Definition 2** (Oscillation). *A solution of (1) is said to be* nonoscillatory *if it is eventually either positive or negative. Otherwise, the solution is called* oscillatory. *In other words, a solution is said to be oscillatory if there exists an increasing divergent sequence $\{\xi_k\}_{k \in \mathbb{N}} \subset [\theta_0, \infty)$ such that $x(\xi_k^+)x(\xi_k^-) \leq 0$ for all $k \in \mathbb{N}$.*

For any given $\varphi \in AC([\rho_t, t], \mathbb{R})$, we define $\|\varphi\| := \sup\{|\varphi(\xi)| : \xi \in [\rho_t, t]\}$.

**Definition 3** (Stability). (i) *The zero solution of (1) is said to be* stable, *if for every $\varepsilon > 0$ and every $\theta \in [\theta_0, \infty)$, there exists $\delta = \delta(\varepsilon, \theta) > 0$ such that any $\varphi \in AC([\rho_\theta, \theta], \mathbb{R})$ with $\|\varphi\| < \delta$ implies $|x(t, \theta, \varphi)| < \varepsilon$ for all $t \in [\theta, \infty)$.*

(ii) *The zero solution of (1) is said to be* uniformly stable, *if for every $\varepsilon > 0$ and every $\theta \in [\theta_0, \infty)$, there exists $\delta = \delta(\varepsilon) > 0$ such that any $\varphi \in AC([\rho_\theta, \theta], \mathbb{R})$ with $\|\varphi\| < \delta$ implies $|x(t, \theta, \varphi)| < \varepsilon$ for all $t \in [\theta, \infty)$.*

*(iii) The zero solution of* (1) *is said to be* asymptotically stable, *if it is stable, and for any $\theta \in [\theta_0, \infty)$ there exists $\delta = \delta(\theta)$ such that any $\varphi \in AC([\rho_\theta, \theta], \mathbb{R})$ with $\|\varphi\| < \delta$ implies $\lim_{t \to \infty} x(t, \theta, \varphi) = 0$.*

The paper is organized as follows: In § 2, we construct the major equipments of the paper which all the results in the sequel will depend on; in § 3, we present our main results, which combine qualitative theory of delay differential equations and qualitative theory of delay differential equations in the absence of retardations in the impulse conditions; in § 4, to conclude the paper, we make our final comments and give a simple example to mention the significance and applicability of the main results. In the sequel, we always assume without furthermore mentioning that $\prod_\emptyset := 1$ and $\sum_\emptyset := 0$.

## 2 Preparatory Results

In this section, we shall introduce several tools required for our main purpose. For simplicity of notation, we let $\theta_{-(\ell+1)} := \varrho_{\theta_0}$. For $i \in \{0, 1, \ldots, \ell\}$, we define

$$\vartheta_k^i := \begin{cases} \theta_0 - \left[\theta_{-(\ell+1)+i+1} - \theta_{-(\ell+1)+i}\right], & k = -1 \\ \theta_0, & k = 0 \\ \vartheta_{k-1}^i + \theta_{k(\ell+1)+i+1} - \theta_{k(\ell+1)+i}, & k \in \mathbb{N}, \end{cases}$$

which explicitly yields

$$\vartheta_k^i = \theta_0 + \sum_{\nu=1}^{k} \left[\theta_{\nu(\ell+1)+i+1} - \theta_{\nu(\ell+1)+i}\right] \quad \text{for } k \in \mathbb{N}.$$

Note that

$$\bigcup_{\nu \in \mathbb{N}_0} [\vartheta_{\nu-1}^i, \vartheta_\nu^i) = \left[\theta_0 - (\theta_{-(\ell+1)+i+1} - \theta_{-(\ell+1)+i}), \infty\right) \quad \text{for each } i \in \{0, 1, \ldots, \ell\}.$$

For $i \in \{0, 1, \ldots, \ell\}$, define

$$\alpha_i : [\vartheta_{-1}^i, \infty) \to \bigcup_{\nu \in \mathbb{N}_0} \left[\theta_{(\nu-1)(\ell+1)+i}, \theta_{(\nu-1)(\ell+1)+i+1}\right)$$

by

$$\alpha_i(t) := \begin{cases} t + [\theta_i - \theta_0] + \sum_{\substack{\theta_0 < \vartheta_j^i \leq t \\ j \in \mathbb{N}}} \left[\theta_{j(\ell+1)+i} - \theta_{(j-1)(\ell+1)+i+1}\right], & t \in [\theta_0, \infty) \\ t - \left[\theta_0 - \theta_{-(\ell+1)+i+1}\right], & t \in [\vartheta_{-1}^i, \theta_0). \end{cases}$$

Then $\alpha_i$ ($i \in \{0, 1, \ldots, \ell\}$) maps the interval $[\vartheta_{k-1}^i, \vartheta_k^i)$ onto $[\theta_{(k-1)(\ell+1)+i}, \theta_{(k-1)(\ell+1)+i+1})$ for each $k \in \mathbb{N}_0$.

Figure 1 is an illustration of the functions $\alpha_i$ ($i = 0, 1, \ldots, \ell$).

In addition to the primary assumptions introduced in the previous section, below, we list two more primary assumptions.

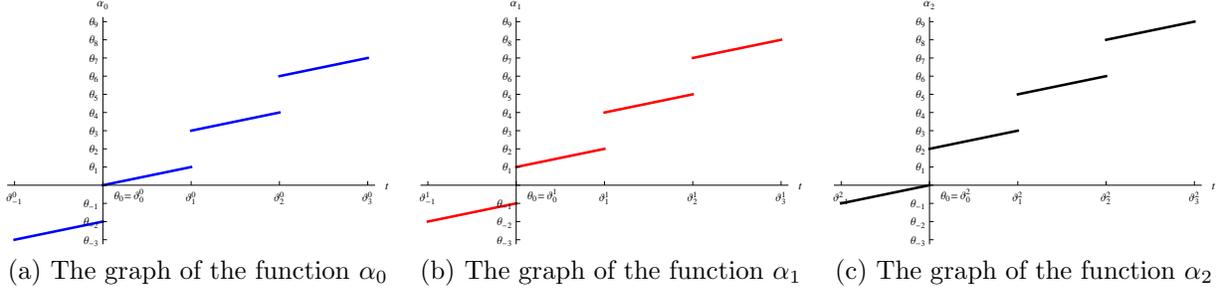

Figure 1: An illustration of the functions $\alpha_i$ $(i = 0, 1, 2)$ with $\ell = 2$

(A5) For each $i \in \{0, 1, \ldots, \ell\}$,
$$t \in \bigcup_{\nu \in \mathbb{N}_0} [\theta_{\nu(\ell+1)+i}, \theta_{\nu(\ell+1)+i+1}) \quad \text{if and only if} \quad \tau(t) \in \bigcup_{\nu \in \mathbb{N}_0} [\theta_{(\nu-1)(\ell+1)+i}, \theta_{(\nu-1)(\ell+1)+i+1}).$$

(A6) There exist functions $\sigma_i : [\theta_0, \infty) \to [\vartheta^i_{-1}, \infty)$ $(i \in \{0, 1, \ldots, \ell\})$ such that
$$\alpha_i(\sigma_i(t)) = \tau(\alpha_i(t)) \quad \text{for all } t \in [\theta_0, \infty).$$

**Remark 1.** Note that for the case $\ell = 0$, we have $\alpha_0(t) = t$ for $t \in [\theta_0, \infty)$, and thus we may let $\sigma_0(t) = \tau(t)$ for $t \in [\theta_0, \infty)$ to have (A5) and (A6) satisfied.

**Lemma 1.** *Assume that (A1)–(A6) hold. Let $x = x(\cdot, \theta_0, \varphi)$ be a solution of (2). Then, for $i \in \{0, 1, \ldots, \ell\}$, the function $y_i : [\theta_0, \infty) \to \mathbb{R}$ defined by*

$$y_i(t) := \left[ \prod_{\substack{\vartheta^i_j \leq t \\ j \in \mathbb{N}_0}} \frac{1}{\lambda_{j(\ell+1)+i}} \right] x(\alpha_i(t)) \quad \text{for } t \in [\theta_0, \infty) \tag{3}$$

*is absolutely continuous on $[\theta_0, \infty)$. Moreover, $i \in \{0, 1, \ldots, \ell\}$, the function $y_i$ is the solution of the initial value problem*

$$\begin{cases} y'_i(t) + \left[ \displaystyle\prod_{\substack{\sigma_i(t) < \vartheta^i_j \leq t \\ j \in \mathbb{N}_0}} \frac{1}{\lambda_{j(\ell+1)+i}} \right] p(\alpha_i(t)) y_i(\sigma_i(t)) = 0 & \text{for } t \in [\theta_0, \infty) \text{ almost everywhere} \\ y_i(t) = \varphi(\alpha_i(t)) & \text{for } t \in [\vartheta^i_{-1}, \theta_0]. \end{cases} \tag{4}$$

*Proof.* Let $i \in \{0, 1, \ldots, \ell\}$, it is easy to see that $y_i$ is absolutely continuous on each of the intervals $[\vartheta^i_{k-1}, \vartheta^i_k)$ for each $k \in \mathbb{N}_0$. Note that

$$\alpha_i(\vartheta^{i+}_k) = \theta^+_{k(\ell+1)+i} \quad \text{and} \quad \alpha_i(\vartheta^{i-}_k) = \theta^-_{(k-1)(\ell+1)+i+1} = \theta^-_{k(\ell+1)+i-\ell} \quad \text{for } k \in \mathbb{N}_0.$$

We now show that $y_i(\vartheta_k^{i+}) = y_i(\vartheta_k^{i-})$ for all $k \in \mathbb{N}$. Indeed, for $k \in \mathbb{N}$, we have

$$y_i(\vartheta_k^i) = \left[\prod_{\substack{\vartheta_j^i \leq \vartheta_k^i \\ j \in \mathbb{N}_0}} \frac{1}{\lambda_{j(\ell+1)+i}}\right] x(\theta_{k(\ell+1)+i}) = \left[\prod_{\substack{\vartheta_j^i \leq \vartheta_k^i \\ j \in \mathbb{N}_0}} \frac{1}{\lambda_{j(\ell+1)+i}}\right] \lambda_{k(\ell+1)+i} x(\theta^-_{k(\ell+1)+i-\ell})$$

$$= \left[\prod_{\substack{\vartheta_j^i < \vartheta_k^i \\ j \in \mathbb{N}_0}} \frac{1}{\lambda_{j(\ell+1)+i}}\right] x(\theta^-_{k(\ell+1)+i-\ell}) = \left[\prod_{\substack{\vartheta_j^i < \vartheta_k^i \\ j \in \mathbb{N}_0}} \frac{1}{\lambda_{j(\ell+1)+i}}\right] x(\theta^-_{(k-1)(\ell+1)+i+1})$$

$$= y_i(\vartheta_k^{i-}),$$

which proves that $y_i(\vartheta_k^{i+}) = y_i(\vartheta_k^{i-})$, and hence $y_i$ is absolutely continuous on $[\vartheta_{-1}^i, \infty)$. On the other hand, we have

$$y_i'(t) + \left[\prod_{\substack{\sigma_i(t) < \vartheta_j^i \leq t \\ j \in \mathbb{N}_0}} \frac{1}{\lambda_{j(\ell+1)+i}}\right] p(\alpha_i(t)) y_i(\sigma_i(t))$$

$$= \left[\prod_{\substack{\vartheta_j^i \leq t \\ j \in \mathbb{N}_0}} \frac{1}{\lambda_{j(\ell+1)+i}}\right] x'(\alpha_i(t)) + \left[\prod_{\substack{\sigma_i(t) < \vartheta_j^i \leq t \\ j \in \mathbb{N}_0}} \frac{1}{\lambda_{j(\ell+1)+i}}\right] p(\alpha_i(t)) \left[\prod_{\substack{\vartheta_j^i \leq \sigma_i(t) \\ j \in \mathbb{N}_0}} \frac{1}{\lambda_{j(\ell+1)+i}}\right] x(\alpha_i(\sigma_i(t)))$$

$$= \left[\prod_{\substack{\vartheta_j^i \leq t \\ j \in \mathbb{N}_0}} \frac{1}{\lambda_{j(\ell+1)+i}}\right] \left\{x'(\alpha_i(t)) + p(\alpha_i(t)) x(\tau(\alpha_i(t)))\right\} = 0$$

for all $t \in [\theta_0, \infty) \setminus \{\theta_k\}_{k \in \mathbb{N}_0}$. Also it is not hard to see that the initial function associated with this equation is $\varphi \circ \alpha_i$ on $[\vartheta_{-1}^i, \theta_0)$. The proof is therefore completed. $\square$

Figure 2 is a graphical illustration of the main idea in the construction of the functions $y_i$ ($i = 0, 1, \ldots, \ell$) from the solution $x$ of (1).

Now, we define $\chi_i$ ($i \in \{0, 1, \ldots, \ell\}$) to be the characteristic function of the interval

$$\bigcup_{\nu \in \mathbb{N}_0} [\theta_{\nu(\ell+1)+i}, \theta_{\nu(\ell+1)+i+1}),$$

i.e.,

$$\chi_i(t) := \begin{cases} 1, & t \in \bigcup_{\nu \in \mathbb{N}_0} [\theta_{\nu(\ell+1)+i}, \theta_{\nu(\ell+1)+i+1}) \\ 0 & \text{otherwise.} \end{cases}$$

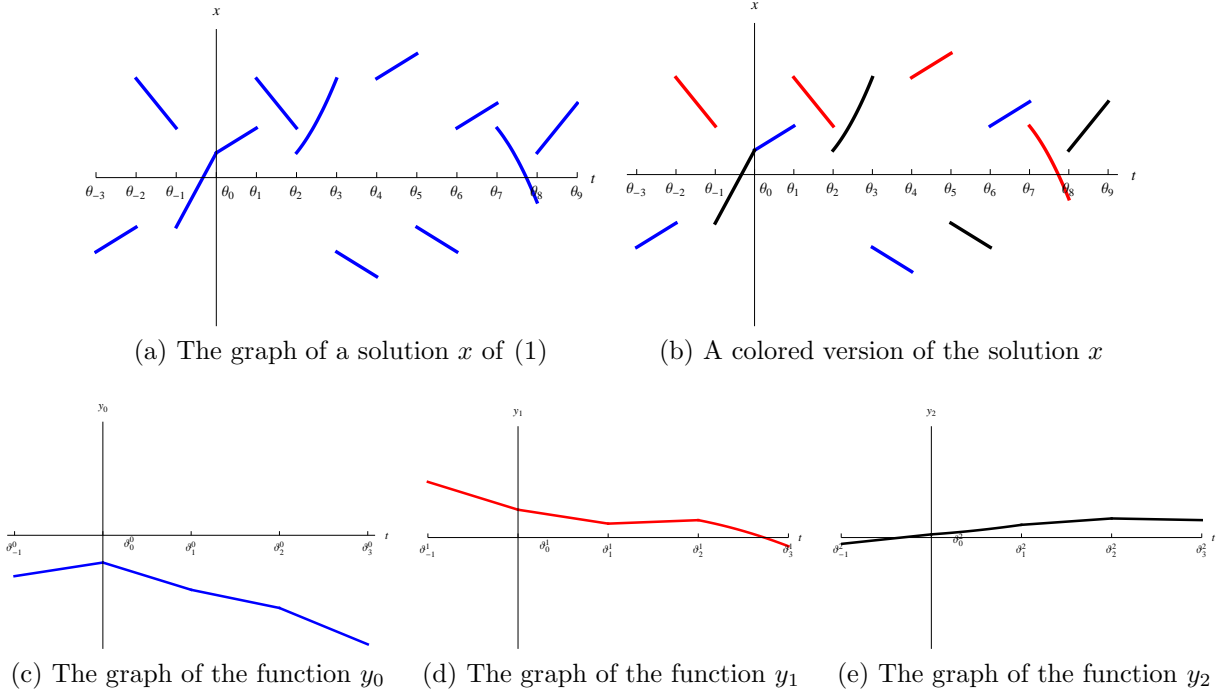

(a) The graph of a solution $x$ of (1)  (b) A colored version of the solution $x$

(c) The graph of the function $y_0$  (d) The graph of the function $y_1$  (e) The graph of the function $y_2$

Figure 2: An illustration of the functions $y_i$ ($i = 0, 1, 2$) with $\ell = 2$

For $i \in \{0, 1, \ldots, \ell\}$, we define the function

$$\beta_i : \bigcup_{\nu \in \mathbb{N}_0} \left[\theta_{(\nu-1)(\ell+1)+i}, \theta_{(\nu-1)(\ell+1)+i+1}\right) \to \left[\vartheta^i_{-1}, \infty\right)$$

by

$$\beta_i(t) := \begin{cases} t - \left\{[\theta_i - \theta_0] + \displaystyle\sum_{\substack{\theta_{j(\ell+1)+i} \leq t \\ j \in \mathbb{N}}} \left[\theta_{j(\ell+1)+i} - \theta_{(j-1)(\ell+1)+i+1}\right]\right\}, & t \in \displaystyle\bigcup_{\nu \in \mathbb{N}_0} \left[\theta_{\nu(\ell+1)+i}, \theta_{\nu(\ell+1)+i+1}\right) \\ t + \left[\theta_0 - \theta_{-(\ell+1)+i+1}\right], & t \in \left[\theta_{-(\ell+1)+i}, \theta_{-(\ell+1)+i+1}\right). \end{cases}$$

It is not hard to see that for each $i \in \{0, 1, \ldots, \ell\}$, $\alpha_i \circ \beta_i$ and $\beta_i \circ \alpha_i$ are the identity mappings on the sets $\bigcup_{\nu \in \mathbb{N}_0} \left[\theta_{(\nu-1)(\ell+1)+i}, \theta_{(\nu-1)(\ell+1)+i+1}\right)$ and $\left[\vartheta^i_{-1}, \infty\right)$, respectively.

**Lemma 2.** *Assume that (A1)–(A4) hold. Let $x$ be a solution of (1) and the functions $y_i : [\theta_0, \infty) \to \mathbb{R}$ ($i \in \{0, 1, \ldots, \ell\}$) be defined by (3). Then*

$$x(t) = \sum_{\mu=0}^{\ell} \chi_\mu(t) \left[\prod_{\substack{\theta_{j(\ell+1)+\mu} \leq t \\ j \in \mathbb{N}_0}} \lambda_{j(\ell+1)+\mu}\right] y_\mu(\beta_\mu(t)) \quad \text{for } t \in [\theta_0, \infty). \tag{5}$$

*Proof.* Let $i \in \{0, 1, \ldots, \ell\}$ and $k \in \mathbb{N}_0$. From (3), we have

$$x\big(\alpha_i(t)\big) = \left[\prod_{\substack{\vartheta_j^i \leq t \\ j \in \mathbb{N}_0}} \lambda_{j(\ell+1)+i}\right] y_i(t) \quad \text{for all } t \in [\theta_0, \infty),$$

which yields

$$x(t) = \left[\prod_{\substack{\vartheta_j^i \leq \beta_i(t) \\ j \in \mathbb{N}_0}} \lambda_{j(\ell+1)+i}\right] y_i\big(\beta_i(t)\big) = \left[\prod_{\substack{\theta_{j(\ell+1)+i} \leq t \\ j \in \mathbb{N}_0}} \lambda_{j(\ell+1)+i}\right] y_i\big(\beta_i(t)\big) \quad \text{for all } t \in [\theta_0, \infty).$$

Therefore, we learn that (5) is true. $\square$

To give an converse analogue of Lemma 1, we need the following additional primary assumption.

(A7) There exist functions $\sigma_i : [\theta_0, \infty) \to [\vartheta_{-1}^i, \infty)$ ($i \in \{0, 1, \ldots, \ell\}$) such that

$$\sigma_i\big(\beta_i(t)\big) = \beta_i\big(\tau(t)\big) \quad \text{for all } t \in [\theta_0, \infty).$$

**Lemma 3.** *Assume that (A1)–(A4), (A5) and (A7) hold. Let $y_i = y_i(\cdot, \theta_0, \varphi_i)$ ($i \in \{0, 1, \ldots, \ell\}$) be solutions of*

$$\begin{cases} y_i'(t) + q_i(t) y_i\big(\sigma_i(t)\big) = 0 & \text{for } t \in [\theta_0, \infty) \\ y_i(t) = \varphi_i(t) & \text{for } t \in [\vartheta_{-1}^i, \theta_0]. \end{cases} \quad (6)$$

*Then, $x$ defined by (5) is a solution of the initial value problem*

$$\begin{cases} x'(t) + \sum_{\mu=0}^{\ell} \chi_\mu(t) \left[\prod_{\substack{\tau(t) < \theta_{j(\ell+1)+\mu} \leq t \\ j \in \mathbb{N}_0}} \lambda_{j(\ell+1)+\mu}\right] q_\mu\big(\beta_\mu(t)\big) x\big(\tau(t)\big) = 0 & \text{for } t \in [\theta_0, \infty) \\ x(\theta_k) = \lambda_k x(\theta_{k-\ell}^-) & \text{for } k \in \mathbb{N}_0 \\ x(t) = \sum_{\mu=0}^{\ell} \chi_\mu(t) \varphi_\mu\big(\beta_\mu(t)\big) & \text{for } t \in [\theta_{-(\ell+1)}, \theta_0). \end{cases} \quad (7)$$

*Proof.* We shall first show that $x$ defined by (5) satisfies the impulse condition in (7). Let $k \in \mathbb{N}_0$,

then we may find $r, i \in \mathbb{N}_0$ such that $k = r(\ell+1) + i$ wand $i \leq \ell$, and thus

$$
\begin{aligned}
x(\theta_k^+) =&\, x\big(\theta_{r(\ell+1)+i}^+\big) = \sum_{\mu=0}^{\ell} \chi_\mu\big(\theta_{r(\ell+1)+i}\big) \left[\prod_{\substack{\theta_{j(\ell+1)+\mu} \leq \theta_{r(\ell+1)+i} \\ j \in \mathbb{N}_0}} \lambda_{j(\ell+1)+\mu}\right] y_\mu\big(\beta_\mu(\theta_{r(\ell+1)+i}^+)\big) \\
=&\, \left[\prod_{\substack{\theta_{j(\ell+1)+i} \leq \theta_{r(\ell+1)+i} \\ j \in \mathbb{N}_0}} \lambda_{j(\ell+1)+i}\right] y_i\big(\beta_i(\theta_{r(\ell+1)+i}^+)\big) = \lambda_{r(\ell+1)+i} \left[\prod_{\substack{\theta_{j(\ell+1)+i} < \theta_{r(\ell+1)+i} \\ j \in \mathbb{N}_0}} \lambda_{j(\ell+1)+i}\right] y_i(\vartheta_r^i) \\
=&\, \lambda_{r(\ell+1)+i} \left[\prod_{\substack{\theta_{j(\ell+1)+i} < \theta_{r(\ell+1)+i} \\ j \in \mathbb{N}_0}} \lambda_{j(\ell+1)+i}\right] y_i\big(\beta_i(\theta_{(r-1)(\ell+1)+i+1}^-)\big) \\
=&\, \lambda_{r(\ell+1)+i} \sum_{\mu=0}^{\ell} \chi_\mu\big(\theta_{(r-1)(\ell+1)+i+1}^-\big) \left[\prod_{\substack{\theta_{j(\ell+1)+\mu} < \theta_{r(\ell+1)+i} \\ j \in \mathbb{N}_0}} \lambda_{j(\ell+1)+\mu}\right] y_\mu\big(\beta_\mu(\theta_{(r-1)(\ell+1)+i+1}^-)\big) \\
=&\, \lambda_{r(\ell+1)+i} x\big(\theta_{(r-1)(\ell+1)+i+1}^-\big) = \lambda_{r(\ell+1)+i} x\big(\theta_{r(\ell+1)+i-\ell}^-\big) = \lambda_k x(\theta_{k-\ell}^-).
\end{aligned}
$$

Now, we show that $x$ defined by (5) satisfies the differential equation condition in (7) but first note that

$$
x'(t) = \sum_{\mu=0}^{\ell} \chi_\mu(t) \left[\prod_{\substack{\theta_{j(\ell+1)+i} \leq t \\ j \in \mathbb{N}_0}} \lambda_{j(\ell+1)+\mu}\right] y'_\mu(\beta_\mu(t))
$$

for all $t \in (\theta_k, \theta_{k+1})$ and all $k \in \mathbb{N}_0$ since the factor of $y_i$ is a step function (its derivative is therefore 0) and $\beta_i$ is a combination of lines of slope 1. Now, we can compute that

$$
\begin{aligned}
& x'(t) + \sum_{\mu=0}^{\ell} \chi_\mu(t) \left[\prod_{\substack{\tau(t) < \theta_{j(\ell+1)+\mu} \leq t \\ j \in \mathbb{N}_0}} \lambda_{j(\ell+1)+\mu}\right] q_\mu\big(\beta_\mu(t)\big) x\big(\tau(t)\big) \\
=& \sum_{\mu=0}^{\ell} \chi_\mu(t) \left[\prod_{\substack{\theta_{j(\ell+1)+\mu} \leq t \\ j \in \mathbb{N}_0}} \lambda_{j(\ell+1)+\mu}\right] y'_\mu(\beta_\mu(t)) + \left\{\sum_{\mu=0}^{\ell} \chi_\mu(t) \left[\prod_{\substack{\tau(t) < \theta_{j(\ell+1)+\mu} \leq t \\ j \in \mathbb{N}_0}} \lambda_{j(\ell+1)+\mu}\right] q_\mu\big(\beta_\mu(t)\big)\right\} \\
& \times \left\{\sum_{\mu=0}^{\ell} \chi_\mu\big(\tau(t)\big) \left[\prod_{\substack{\theta_{j(\ell+1)+\mu} \leq \tau(t) \\ j \in \mathbb{N}_0}} \lambda_{j(\ell+1)+\mu}\right] y_\mu\big(\beta_\mu(\tau(t))\big)\right\}
\end{aligned}
$$

$$= \sum_{\mu=0}^{\ell} \chi_\mu(t) \left[ \prod_{\substack{\theta_{j(\ell+1)+\mu} \leq t \\ j \in \mathbb{N}_0}} \lambda_{j(\ell+1)+\mu} \right] y'_\mu(\beta_\mu(t)) + \left\{ \sum_{\mu=0}^{\ell} \chi_\mu(t) \left[ \prod_{\substack{\tau(t) < \theta_{j(\ell+1)+\mu} \leq t \\ j \in \mathbb{N}_0}} \lambda_{j(\ell+1)+\mu} \right] q_\mu(\beta_\mu(t)) \right\}$$

$$\times \left\{ \sum_{\mu=0}^{\ell} \chi_\mu(t) \left[ \prod_{\substack{\theta_{j(\ell+1)+\mu} \leq \tau(t) \\ j \in \mathbb{N}_0}} \lambda_{j(\ell+1)+\mu} \right] y_\mu(\sigma_\mu(\beta_\mu(t))) \right\}$$

$$= \sum_{\mu=0}^{\ell} \chi_\mu(t) \left[ \prod_{\substack{\theta_{j(\ell+1)+\mu} \leq t \\ j \in \mathbb{N}_0}} \lambda_{j(\ell+1)+\mu} \right] \left\{ y'_\mu(\beta_\mu(t)) + q_\mu(\beta_\mu(t)) y_\mu(\sigma_\mu(\beta_\mu(t))) \right\} = 0,$$

where we have used the condition (A5) and (A7) in the third step above, and the property of the characteristic function $\chi_i$ in the last step for collecting the terms in a single sum. It is not hard to see that the solution $x$ initiates with the function $\sum_{\mu=0}^{\ell} \chi_\mu \times (\varphi_\mu \circ \beta_\mu)$ on $[\theta_{-(\ell+1)}, \theta_0)$. The proof is therefore completed. □

**Remark 2.** Assuming that (A1)–(A6) and (A7) are true, if we replace $q_i$ and $\varphi_i$ in Lemma 3 with $\left[ \prod_{\sigma_i(t) < \vartheta^i_j \leq t, j \in \mathbb{N}_0} (1/\lambda_{j(\ell+1)+i}) \right] \times (p \circ \alpha_i)$ and $\varphi \circ \alpha_i$, respectively, then we obtain the converse of Lemma 1.

## 3 Main Results

In this section, by combining the conclusions of Lemmas 1–3 and Remark 2, we establish the following fundamental theorem, which builds a bridge between impulsive differential equations and nonimpulsive differential equations.

**Theorem 1.** *Assume that (A1)–(A7) hold. Then the following assertions are true:*

(i) *If $x(\cdot, \theta_0, \varphi)$ is a solution of (1), then $y_i(\cdot, \theta_0, \varphi \circ \alpha_i)$ ($i \in \{0, 1, \ldots, \ell\}$) defined by (3) is a solution of (4).*

(ii) *If $y_i(\cdot, \theta_0, \varphi \circ \alpha_i)$ ($i \in \{0, 1, \ldots, \ell\}$) are solutions of (4), then $x(\cdot, \theta_0, \varphi)$ defined by (5) is a solution of (1).*

*Proof.* The proof of the parts (i) and (ii) follow from Lemmas 1–3 and Remark 2, respectively. □

Next, we introduce an additional assumption to study oscillation of solutions to (1).

(A8) $\{\lambda_k\}_{k \in \mathbb{N}_0}$ is a sequence of positive reals.

**Theorem 2.** *Assume that (A1)–(A3), (A5)–(A8) hold. Then every solution of (1) is oscillatory on $\bigcup_{\nu \in \mathbb{N}_0} [\theta_{\nu(\ell+1)+i}, \theta_{\nu(\ell+1)+i+1})$ for some $i \in \{0, 1, \ldots, \ell\}$ if and only if every solution of (4) is oscillatory for the same $i$.*

*Proof.* Suppose that $x$ is a solution of (1), which is oscillatory on $\bigcup_{\nu \in \mathbb{N}_0} \left[\theta_{\nu(\ell+1)+i}, \theta_{\nu(\ell+1)+i+1}\right)$ for some $i \in \{0, 1, \ldots, \ell\}$. Then for the same $i$, the function $y_i$ defined by (3) is oscillatory on $[\theta_0, \infty)$ since the condition (A8) implies that the transform is oscillation invariant and the function $\alpha_i$ maps the union of intervals $\bigcup_{\nu \in \mathbb{N}_0} \left[\theta_{\nu(\ell+1)+i}, \theta_{\nu(\ell+1)+i+1}\right)$ onto the half-line $[\theta_0, \infty)$. Clearly, $y_i$ is a solution of (4) for the same $i$. Conversely, if for some $i \in \{0, 1, \ldots, \ell\}$, $y_i$ is an oscillating solution of (4), then $x$ defined by (5) is oscillatory on $\bigcup_{\nu \in \mathbb{N}_0} \left[\theta_{\nu(\ell+1)+i}, \theta_{\nu(\ell+1)+i+1}\right)$ and is a solution to (1) because of the fact that $\beta_i$ maps the half-line $[\theta_0, \infty)$ onto the union of intervals $\bigcup_{\nu \in \mathbb{N}_0} \left[\theta_{\nu(\ell+1)+i}, \theta_{\nu(\ell+1)+i+1}\right)$. The proof is therefore completed. □

**Corollary 1.** Assume that (A1)–(A3), (A5)–(A8) hold. Then every solution of (1) is oscillatory if every solution of (4) is oscillatory some $i \in \{0, 1, \ldots, \ell\}$.

**Remark 3.** If every solution of (4) is nonoscillatory for all $i \in \{0, 1, \ldots, \ell\}$, then this does not mean that (1) is nonoscillatory too. For instance if the solutions $y_i$ and $y_j$ of (4) are of different signs, then the solution $x$ of (1), which is defined by (5), has different signs on the union intervals $\bigcup_{\nu \in \mathbb{N}_0} \left[\theta_{\nu(\ell+1)+i}, \theta_{\nu(\ell+1)+i+1}\right)$ and $\bigcup_{\nu \in \mathbb{N}_0} \left[\theta_{\nu(\ell+1)+j}, \theta_{\nu(\ell+1)+j+1}\right)$, and therefore, $x$ is oscillatory.

**Remark 4.** If the sequence $\{\lambda_k\}_{k \in \mathbb{N}_0}$ has a negative subsequence $\{\lambda_{k_l}\}_{l \in \mathbb{N}_0}$, then every solution of (1) changes sign at the impulse points $\{\theta_{k_l}\}_{l \in \mathbb{N}_0}$, i.e., every solution of (1) oscillates.

Next to study stability of (1), we give two more assumptions.

(A9) For any $s \geq \theta_0$, there exists a positive constant $M = M(s)$ such that

$$\sum_{k=0}^{\ell} \left| \prod_{\substack{s < \vartheta_j^k \leq t \\ j \in \mathbb{N}_0}} \frac{1}{\lambda_{j(\ell+1)+k}} \right| \leq M \quad \text{for all } t \geq s.$$

(A10) There exists a positive constant $M$ such that

$$\sum_{k=0}^{\ell} \left| \prod_{\substack{s < \vartheta_j^k \leq t \\ j \in \mathbb{N}_0}} \frac{1}{\lambda_{j(\ell+1)+k}} \right| \leq M \quad \text{for all } s, t \geq \theta_0 \text{ with } t \geq s.$$

**Theorem 3.** *Assume that (A1)–(A7) hold. Then the following assertions are true:*

(i) *Suppose that (A9) holds. If the zero solution of (4) is stable for all $i \in \{0, 1, \cdots, \ell\}$, then so does the zero solution of (1).*

(ii) *Suppose that (A10) holds. If the zero solution of (4) is uniformly stable for all $i \in \{0, 1, \cdots, \ell\}$, then so does the zero solution of (1).*

(iii) *Suppose that (A9) holds. If the zero solution of (4) is asymptotically stable for all $i \in \{0, 1, \cdots, \ell\}$, then so does the zero solution of (1).*

*Proof.* We shall only give a proof for the part (i) since the proof of the parts (ii) and (iii) follow similar arguments. Let $\varepsilon > 0$ and $s \in [\theta_0, \infty)$. From the hypothesis, for $\varepsilon_i > 0$ with $\max_{i \in \{0,1,\ldots,\ell\}} \{\varepsilon_i\} \leq \varepsilon$ and $s_i \in [\theta_0, \infty)$ with $\max_{i \in \{0,1,\ldots,\ell\}} \{\alpha_i(s_i)\} \leq s$, we may find $\delta_i = \delta_i(\varepsilon_i, s_i) > 0$ such that $\varphi_i \in AC([\rho^i_{s_i}, s_i], \mathbb{R})$, where $\rho^i_t := \inf\{\sigma_i(\xi) : \xi \in [t, \infty)\}$ for $t \in [\theta_0, \infty)$, with $\|\varphi_i\| < \delta_i$, which implies

$$|y_i(t)| < \frac{\varepsilon_i}{M(\ell+1)} \quad \text{for all } t \in [s_i, \infty).$$

From (5) and (A9), we have

$$|x(t)| = \left| \sum_{\mu=0}^{\ell} \chi_\mu(t) \left[ \prod_{\substack{s < \theta_{j(\ell+1)+\mu} \leq t \\ j \in \mathbb{N}_0}} \lambda_{j(\ell+1)+\mu} \right] y_\mu(\beta_\mu(t)) \right|$$

$$= \left| \sum_{\mu=0}^{\ell} \chi_\mu(t) \left[ \prod_{\substack{\beta_\mu(s) < \vartheta^\mu_j \leq \beta_\mu(t) \\ j \in \mathbb{N}_0}} \lambda_{j(\ell+1)+\mu} \right] y_\mu(\beta_\mu(t)) \right|$$

$$\leq \sum_{\mu=0}^{\ell} \chi_\mu(t) \left[ \prod_{\substack{\beta_\mu(s) < \vartheta^\mu_j \leq \beta_\mu(t) \\ j \in \mathbb{N}_0}} \lambda_{j(\ell+1)+\mu} \right] |y_\mu(\beta_\mu(t))|$$

$$\leq \left( \sum_{\mu=0}^{\ell} \left[ \prod_{\substack{\beta_\mu(s) < \vartheta^\mu_j \leq \beta_\mu(t) \\ j \in \mathbb{N}_0}} \lambda_{j(\ell+1)+\mu} \right] \right) \left( \sum_{\mu=0}^{\ell} |y_\mu(\beta_\mu(t))| \right)$$

$$< M \sum_{\mu=0}^{\ell} \frac{\varepsilon_i}{(\ell+1)M} \leq \varepsilon$$

for all $t \in [s, \infty)$. The zero solution of (1) is therefore stable. $\square$

(A11) For any $s \geq \theta_0$, there exists a positive constant $N = N(s)$ such that

$$\sum_{k=0}^{\ell} \left| \prod_{\substack{s < \theta_{j(\ell+1)+k} \leq t \\ j \in \mathbb{N}_0}} \lambda_{j(\ell+1)+k} \right| \leq N \quad \text{for all } t \geq s.$$

(A12) There exists a positive constant $N$ such that

$$\sum_{k=0}^{\ell} \left| \prod_{\substack{s < \theta_{j(\ell+1)+k} \leq t \\ j \in \mathbb{N}_0}} \lambda_{j(\ell+1)+k} \right| \leq N \quad \text{for all } s, t \geq \theta_0 \text{ with } t \geq s.$$

As a dual version of Theorem 3, we give the following result. The proof is omitted since it makes use of almost the same arguments.

**Theorem 4.** *Assume that (A1)–(A7) hold. Then the following assertions are true:*

(i) *Suppose that (A11) holds. If the zero solution of (1) is stable, then so does the zero solution of (4) for all $i \in \{0, 1, \ldots, \ell\}$.*

(ii) *Suppose that (A12) holds. If the zero solution of (1) is uniformly stable, then so does the zero solution of (4) for all $i \in \{0, 1, \ldots, \ell\}$.*

(iii) *Suppose that (A11) holds. If the zero solution of (1) is asymptotically stable, then so does the zero solution of (4) for all $i \in \{0, 1, \ldots, \ell\}$.*

We conclude this section with the following corollary which combines the conclusions of Theorem 3 and Theorem 4.

**Corollary 2.** Assume that (A1)–(A7) hold. Then the following assertions are true:

(i) Suppose that (A9) and (A11) hold. The zero solution of (1) is stable if and only if so does the zero solution of (4) for all $i \in \{0, 1, \ldots, \ell\}$.

(ii) Suppose that (A10) and (A12) holds. If the zero solution of (1) is uniformly stable if and only if so does the zero solution of (4) for all $i \in \{0, 1, \ldots, \ell\}$.

(iii) Suppose that (A9) and (A11) holds. If the zero solution of (1) is asymptotically stable if and only if so does the zero solution of (4) for all $i \in \{0, 1, \ldots, \ell\}$.

# 4 Final Comments

In this section, we shall make certain comments about how to extend the method, and then present a simple example to illustrate the applicability of our results.

The method introduced here allows the study of qualitative behaviour of equations of the form

$$\begin{cases} x'(t) + \sum_{i=1}^{n} p_i(t) x\big(\tau_i(t)\big) = 0 & \text{for } t \in [\theta_0, \infty) \setminus \{\theta_k\}_{k \in \mathbb{N}_0} \\ x(\theta_k) = \lambda_k x(\theta_{k-\ell}^-) & \text{for } k \in \mathbb{N}_0, \end{cases}$$

however just for simplicity of notation, we have restricted our attention to (1). On the other hand, one can obtain straightforwardly the corresponding results for impulsive differential equations with left-continuous solutions, i.e., $x(\theta_k) = x(\theta_k^-)$ for all $k \in \mathbb{N}_0$. Assuming that the coefficient $p$ is nonnegative, one would wish to have results for left-continuous solutions since eventually positive solutions become nonincreasing in this case. Since our results introduced in § 3 do not require any sign condition on the coefficient $p$, we preferred right-continuous solutions.

Suppose there exist $n \in \mathbb{N}$ and a set of additional points $\{\theta_{-n(\ell+1)}, \theta_{-n(\ell+1)+1}, \ldots, \theta_{-n(\ell+1)}\} \subset [\rho_{\theta_0}, \theta_{-\ell}]$ with the convention $\theta_{-n(\ell+1)} := \rho_{\theta_0}$. Then, we may let

$$\vartheta_k^i := \begin{cases} \vartheta_{k+1}^i - [\theta_{k(\ell+1)+i+1} - \theta_{k(\ell+1)+i}], & k \in \{-n, -n+1, \ldots, -1\} \\ \theta_0, & k = 0 \\ \vartheta_{k-1}^i + \theta_{k(\ell+1)+i+1} - \theta_{k(\ell+1)+i}, & k \in \mathbb{N}, \end{cases}$$

and extend the domain of the function $\alpha_i$ in a very natural way such that

$$\alpha_i(t) := \begin{cases} t + [\theta_i - \theta_0] + \displaystyle\sum_{\substack{\theta_0 < \vartheta_j^i \leq t \\ j \in \mathbb{N}}} [\theta_{j(\ell+1)+i} - \theta_{(j-1)(\ell+1)+i+1}], & t \in [\theta_0, \infty) \\[2em] t - \left[\theta_0 - \theta_{-(j+1)(\ell+1)+i+1}\right]_{\substack{\vartheta_{-(j+1)}^i < t \leq \vartheta_{-j}^i \\ j \in \mathbb{N}}} + \displaystyle\sum_{\substack{t < \vartheta_{-j}^i \leq \vartheta_{-1}^i \\ j \in \mathbb{N}}} [\theta_{-j(\ell+1)+i+1} - \theta_{-j(\ell+1)+i}], & t \in [\vartheta_{-n}^i, \theta_0) \end{cases}$$

for $t \in [\vartheta_{-n}^i, \infty)$, and similarly we can extend the function

$$\beta_i : \bigcup_{\nu \in \mathbb{N}_0} [\theta_{(\nu-n)(\ell+1)+i}, \theta_{(\nu-n)(\ell+1)+i+1}) \to [\vartheta_{-n}^i, \infty).$$

by

$$\beta_i(t) := \begin{cases} t - \left\{[\theta_i - \theta_0] + \displaystyle\sum_{\substack{\theta_{j(\ell+1)+i} \leq t \\ j \in \mathbb{N}}} [\theta_{j(\ell+1)+i} - \theta_{(j-1)(\ell+1)+i+1}]\right\}, & t \in \bigcup_{\nu \in \mathbb{N}_0} [\theta_{\nu(\ell+1)+i}, \theta_{\nu(\ell+1)+i+1}) \\[2em] t + \left\{\left[\theta_0 - \theta_{-(j+1)(\ell+1)+i+1}\right]_{\substack{\theta_{-(j+1)(\ell+1)+i} < t \leq \theta_{-j(\ell+1)+i} \\ j \in \mathbb{N}}} \right. \\ \qquad \left. - \displaystyle\sum_{\substack{t < \theta_{-j(\ell+1)+i} \leq \theta_{-(\ell+1)+i} \\ j \in \mathbb{N}}} [\theta_{-j(\ell+1)+i+1} - \theta_{-j(\ell+1)+i}]\right\}, & t \in \bigcup_{\nu=1}^{n} [\theta_{-\nu(\ell+1)+i}, \theta_{-\nu(\ell+1)+i+1}). \end{cases}$$

Therefore our results in § 3 are still true if we replace (A5)–(A7) with the following assumptions, respectively:

(A13) For each $i \in \{0, 1, \ldots, \ell\}$,

$$t \in \bigcup_{\nu \in \mathbb{N}_0} [\theta_{\nu(\ell+1)+i}, \theta_{\nu(\ell+1)+i+1}) \quad \text{if and only if} \quad \tau(t) \in \bigcup_{\nu \in \mathbb{N}_0} [\theta_{(\nu-n)(\ell+1)+i}, \theta_{(\nu-n)(\ell+1)+i+1}).$$

(A14) There exist functions $\sigma_i : [\theta_0, \infty) \to [\vartheta_{-n}^i, \infty)$ such that

$$\alpha_i(\sigma_i(t)) = \tau(\alpha_i(t)) \quad \text{for all } t \in [\theta_0, \infty)$$

and for $i \in \{0, 1, \ldots, \ell\}$.

(A15) There exist functions $\sigma_i : [\theta_0, \infty) \to [\vartheta^i_{-n}, \infty)$ such that
$$\sigma_i\bigl(\beta_i(t)\bigr) = \beta_i\bigl(\tau(t)\bigr) \quad \text{for all } t \in [\theta_0, \infty)$$
and for $i \in \{0, 1, \ldots, \ell\}$.

The following example considers the extended versions of the results in § 3.

**Example 1.** *Consider*
$$\begin{cases} x'(t) + px\bigl(t - n(\ell+1)\bigr) = 0 & \text{for } t \in \mathbb{R}_0^+ \setminus \mathbb{N}_0 \\ x(k) = \lambda_{\mathrm{mod}(k,\ell+1)} x\bigl((k-\ell)^-\bigr) & \text{for } k \in \mathbb{N}_0, \end{cases} \qquad (8)$$

*where $n \in \mathbb{N}$, $\ell \in \mathbb{N}_0$, $p \in \mathbb{R}$ and $\lambda_i \in \mathbb{R}\setminus\{0\}$ for $i \in \{0, 1, \ldots, \ell\}$. For this equation, we have $\tau(t) = t - n(\ell+1)$ for $t \in \mathbb{R}_0^+$, which satisfies the condition (A13) if we let $\theta_k = k$ for $k \in \{-n(\ell+1), -n(\ell+1)+1, \ldots\}$. For $i \in \{0, 1, \ldots, \ell\}$, we can easily compute that $\vartheta^i_k = k$ for $k \in \{-n, -n+1, \ldots\}$, and for $t \in [-n, \infty)$, we have*

$$\alpha_i(t) = \begin{cases} t + i + \displaystyle\sum_{\substack{j \leq t \\ j \in \mathbb{N}}} \ell, & t \in [0, \infty) \\[2ex] t + \left(\displaystyle\sum_{\substack{-(j+1) < t \leq -j \\ j \in \mathbb{N}}} -(j+1)(\ell+1) + i + 1\right) + \displaystyle\sum_{\substack{t < -j \leq -1 \\ j \in \mathbb{N}}} 1, & t \in [-n, 0) \end{cases}$$
$$= t + i + \lfloor t \rfloor \ell,$$

*where $\lfloor \cdot \rfloor$ denotes the least integer function. For $i \in \{0, 1, \ldots, \ell\}$, letting $\sigma_i(t) = t - n$ for $t \in [-n, \infty)$, we see that*

$$\begin{aligned} \alpha_i\bigl(\sigma_i(t)\bigr) &= \sigma_i(t) + i + \lfloor \sigma_i(t) \rfloor \ell = (t-n) + i + \lfloor (t-n) \rfloor \ell \\ &= (t-n) + i + \bigl(\lfloor t \rfloor - n\bigr)\ell = t + i + \lfloor t \rfloor \ell - n(\ell+1) \\ &= \tau\bigl(t + i + \lfloor t \rfloor \ell\bigr) = \tau(\sigma_i(t)) \end{aligned}$$

*for all $t \in [-n, \infty)$, which shows that (A14) is satisfied. One can show very similarly that (A15) is also satisfied. Then, from Lemma 1, for $i \in \{0, 1, \ldots, \ell\}$, the associated differential equations have the form*
$$y'_i(t) + \frac{p}{\lambda_i^n} y_i(t-n) = 0 \quad \text{for } t \in [0, \infty) \text{ almost everywhere}$$
*since*
$$\prod_{\substack{t-n < j \leq t \\ j \in \mathbb{N}_0}} \frac{1}{\lambda_{\mathrm{mod}(j(\ell+1)+i,\ell+1)}} = \frac{1}{\lambda_i^n} \quad \text{for all } t \in [0, \infty).$$

*Clearly, (A9)–(A12) are satisfied. Therefore, by the extended version of Theorem 2 (see also Corollary 1), every solution of (8) is oscillatory if*
$$\lambda_i > 0 \quad \text{for all } i \in \{0, 1, \ldots, n\} \quad \text{and} \quad n\frac{p}{\lambda_i^n} > \frac{1}{\mathrm{e}} \quad \text{for some } i \in \{0, 1, \ldots, n\}$$

(see [17]), and by the extended version of Theorem 3 (see also Corollary 2), every solution tends to zero asymptotically if

$$\frac{p}{\lambda_i^n} > 0 \quad and \quad n\frac{p}{\lambda_i^n} \leq \frac{\pi}{2} \quad for\ all\ i \in \{0, 1, \ldots, n\}$$

(see [13, 18]).